\documentstyle[preprint,aps]{revtex}
\widetext
\draft
\tighten

\begin{document}

\title{On an ambiguity in the concept of partial and total derivatives in
classical analysis}

\bigskip

\author{{\bf Andrew E. Chubykalo,   Rolando A.  Flores\thanks{Centro de
Estudios Multidisciplinarios,  Universidad Aut\'onoma de Zacatecas } and
Juan A.  P\'erez\thanks{Centro Regional de Estudios Nucleares, Zacatecas}}
}

\address {Escuela de F\'{\i}sica, Universidad Aut\'onoma de Zacatecas \\
Apartado Postal C-580\, Zacatecas 98068, ZAC., M\'exico}

\date{\today}

\maketitle


\baselineskip 7mm
\bigskip
\bigskip
\bigskip
\bigskip


\begin{abstract}
Ambiguity is shown in the context
of the differential calculus of several variables and
with the
help of the language of category theory, a way to solve it in its most
general form is offered. It is also shown that this new definition is
related to other well-known definitions in the literature.
\end{abstract}

\section{Introduction}

The difference between the functions:
$$
E [x_1(t),\ldots, x_{n-1}(t), t]  = \;_{{\rm def}}E[{\bf r}(t),t],\quad
E(x_1,\ldots,x_{n-1},t) = \;_{{\rm def}}E({\bf r}, t)
$$
is usually  not remarked in the literature, and for this reason we can
often write down meaningless symbols like:
\begin{equation}
\frac{\partial}{\partial t}E[{\bf r}(t), t],
\end{equation}
and
\begin{equation}
\frac{d}{dt}E({\bf r},t).
\end{equation}

Ambiguities in the ``notation" for partial differentiation has been
remarked by Arnold [1] p. 226 (p. 258 in English translation) without
further development. The  symbols (1), (2) are meaningless, because the
process denoted by the operator of {\it partial}  differentiation can be
applied only to functions of several {\it independent} variables and  $E
[{\bf r}(t), t]$ is not {\it such} a function.  Meanwhile, the operator of
{\it total} differentiation with respect to given variable can  be
formally applied to functions of one variable only.  However,  we have a
well-known formula to relate both concepts:  \begin{equation}
\frac{d}{dt}E  = ({\bf V}\cdot\nabla)E + \frac{\partial}{\partial t}E
\end{equation} (here ${\bf V} = \frac{d{\bf r}}{dt}$).

Let us show that, in this form, Eq. (3) cannot be correct. What  is the
correct argument for the symbol $E$ in both sides? If we say that the
correct argument for both sides is $[{\bf r}(t), t]$ we get the chain of
symbols  (1), but in this case, the operator of a partial
differentiation would indicate that we must construct a new function in
the form $(\partial E/\partial t)$, hence we use the following procedure:

\begin{equation}
\lim\limits_{\Delta t\rightarrow 0}\left\{\frac{E\left[{\bf r}(t) + \Delta
t\frac{d{\bf r}(t)}{dt},\; t+\Delta t\right] - E[{\bf r}(t), t]}{\Delta
t}\right\}.
\end{equation}
But this  is the definition of total
differentiation! Thus, the symbols of  total and of partial
differentiation denote the same process, therefore,  because $E$ is the
same function on both sides of the equation, we get:
\begin{equation}
({\bf V}\cdot\nabla)E[{\bf r}(t), t]  = 0
\end{equation}
always.  But even if the procedure which we followed were correct
(which it is not, of course!),  this equation is not correct for $E$ as a
function of the functions ${\bf r}(t)$, because the partial
differentiation would involve increments of the functions ${\bf r}(t)$ in
the form ${\bf r}(t) + \Delta{\bf r}(t)$ and we do not know how we must
interpret this increment because we have two options: {\it either} $\Delta
{\bf r}(t) = {\bf r}(t) - {\bf r}^*(t)$, {\it or} $\Delta{\bf r}(t) =
{\bf r}(t) - {\bf r}(t^*)$.  Both are different processes because the
first one involves changes in the functional form of the functions ${\bf
r}(t)$, while the second involves changes in the position along the path
defined by ${\bf r} = {\bf r}(t)$ but preserving the same functional form.
Hence, it is clear that we have here different concepts.  If we remember
the definition of partial differentiation, we can see where the mistake is:
``{\it the symbol:  $\frac{\partial}{\partial t}E({\bf r}, t)$ means that
we take the variations of $t$ when the values of} {\bf r} {\it are
constant}". It means that we make the only change $t + \Delta t$ in the
function.  But this is only possible if the coordinates {\bf r} are
independent from $t$.  Hence, we can see that the correct argument cannot
be $[{\bf r}(t), t]$, because, as we have shown, this supposition leads to
the incorrect result (5). If we make the other supposition, that the
correct argument is $({\bf r}, t)$ we can get the same conclusion, i.e.,
equation (5).  Hence, {\it none of these suppositions is correct}.  What
is the solution, then?  Actually, in the
equation (3) we have two {\it different} functions: on the left hand side
we have the function $E[{\bf r}(t),t]$ defined on a {\it curve} in a
$n$-surface and on the right hand side we have the function $E({\bf r},
t)$ defined on the {\it all}  $n$-surface, which obviously are {\it
quite} different functions, while  we have a limiting procedure to get a
unification of concepts in the realm of functions of one variable.

Now let us introduce the following notation:
\begin{equation}
f=E\circ {\bf p},
\end{equation}
where the symbol ``$\circ$" means a composition of functions and where
$$
E:R^n\rightarrow R, \quad {\bf p}:R\rightarrow R^{n},\quad f:R\rightarrow
R.  $$
It is clear that ${\bf p}={\bf p}(t)=\{x_1(t),\ldots,x_{n-1}(t),t\}=\{{\bf
r}(t),t\}$ is a curve which lies on the $n$-surface where the function $E$
is defined.

Hence we can write down the equation:
$$ \frac{d}{dt}f=\lim\limits_{x_i\rightarrow x_i(t)}\left\{({\bf
V}\cdot\nabla)E + \frac{\partial E}{\partial t}\right\}
$$
which shows our point more clearly: {\it the functions in both sides ($f$
and $E$) are different functions}. Of course, we suppose that the
components of the vector {\bf V} tend to derivatives $\frac{dx_i}{dt}$ in
the limit.  But here is where our grammatical distinction appears:  the
right hand side is evaluated in all points along the curve ${\bf p}(t)$,
that is:  $$ ({\bf V}\cdot\nabla)E\Bigl|_{x_i=x_i(t)}+ \frac{\partial
E}{\partial t}\Bigl|_{x_i=x_i(t)}.  $$

Let us explain the distinction as follows: the operator of the total
differentiation is  just  a differentiation of a function which can
depend on one independent variable, and the operator of the partial
differentiation is just a partial differentiation of a function which can
depend on several independent variables.
An obvious question immediately  arises:  {\it what is the relation
between these domains}?  Obviously, the function of one variable is one
entity and the function of several variables is a different one. The
relation lies in the evaluation of the function obtained by partial
differentiation in points {\it along} the curve ${\bf p}(t)$. Or in more
general terms, we must have the validity of the following condition: for
all $\varepsilon >0$ there is a $\delta[\varepsilon, {\bf p}(t)]>0$ such
that if we take a point in the ball:  $$ |\hat{\bf r}-{\bf
p}(t)|<\delta[\varepsilon,{\bf p}(t)], $$ where $\hat{\bf
r}=\{x_1,\ldots,x_{n-1},t\}=\{{\bf r},t\}$, then $$
\left|\left[\Bigl\{({\bf V}\cdot\nabla)E\Bigr\}({\bf
r},t)+\left\{\frac{\partial E}{\partial t}\right\}({\bf r},t)\right]-
\left[\Bigl\{({\bf V}\cdot\nabla)E\Bigr\}[{\bf
p}(t)]+\left\{\frac{\partial E}{\partial t}\right\}[{\bf
 p}(t)]\right]\right|<\varepsilon, $$
where, of course, $E=E(\hat{\bf r})=E({\bf r},t)$.

We have not supposed, of course, that we have an uniform continuity.
The abbreviated form of this condition is:
\begin{equation}
\frac{d}{dt}f(t) =
\lim\limits_{\hat{\bf r}\rightarrow {\bf p}(t)}\left\{({\bf V}\cdot
\nabla)E({\bf r}, t) + \frac{\partial}{\partial t}E({\bf r},t)\right\}.
\end{equation}

The distinction  between\footnote{or, that
is the same, between $E(\hat{\bf r})$ and $E[{\bf p}(t)]$} $E({\bf r}, t)$
and $E[{\bf r}(t),t]$ is important in some physical contexts, as it is
shown in [2] (see, especially, Eq.  (28)).  The grammatical distinction
is that the realm of functions of one independent variable is not the same
as the realm of functions of several independent variables, and that the
relation between these two realms is given by a limitation procedure.

\section{Some remarks related to the functional equation (6)}

What conditions must  the relation $f=E\circ{\bf p}$ satisfy to make
sense? It is obviously that all its elements $f$, $E$, and $p$ have to
exist, and we, in fact, must write down the more general relation:
\begin{equation}
f(t)=\lim\limits_{\hat{\bf r}\rightarrow{\bf p}}E({\bf r},t).
\end{equation}
It means that the function $E$ must be continuous in all points of the
curve {\bf p}.

We have to consider seven cases\footnote{When we know all three functions,
we must only check that the relation (7) is valid. This is trivial.}:
\begin{quotation}

1. Two functional form $E$ and {\bf p} are known: This is the $\{E,{\bf
p}\}$-case;

2. $E$ and $f$ are known: $\{E,f\}$-case;

3. {\bf p} and $f$ are known: $\{{\bf p},f\}$-case;

4. Only $E$ is known: $\{E\}$-case;

5. Only $f$ is known: $\{f\}$-case;

6. Only {\bf p} is known: $\{{\bf p}\}$-case;

7. All function are unknown: $\{\}$-case.

\end{quotation}

In the $\{E,{\bf p}\}$-, $\{E,f\}$- and $\{{\bf p},f\}$-cases we can
define one of the functions in terms of the other two functions, for
example, in $\{E,f\}$-case we define ${\bf p}$ etc. In $\{E\}$-, $\{f\}$-
and $\{{\bf p}\}$-cases one can show that it is possible to define the
other functions under certain conditions. Let us make a brief review of
these classes.

\begin{quotation}
\underline{$\{E\}$-case:} In this case, we only know the form of $E$, and
we need to define the forms of the other two functions. We suppose
that:  $E\in C^1(R^n,R)$, ${\bf p}\in C^1(R,R^n)$, $f\in C^1(R,R)$. Now we
write down our defining equation in the form
\begin{equation}
\frac{df}{dt}=\lim\limits_{\hat{\bf r}\rightarrow{\bf
p}}\left\{\sum\limits_{i=1}^{n-1}V_i(\hat{\bf r})\frac{\partial
E}{\partial x_i}+ \frac{\partial E}{\partial t}\right\},
\end{equation}
and we propose the following two equations:

\begin{equation}
({\rm a})\;\frac{df}{dt}=
\lim\limits_{\hat{\bf r}\rightarrow{\bf
p}}\frac{\partial E}{\partial t}\qquad
{\rm and}\qquad({\rm b})\;
V_i(\hat{\bf r})=\sum\limits_{j=1}^{n-1}b_{ij}\frac{\partial E}{\partial
x_j},
\end{equation}
where $b_{ij}$ is a skew-symmetrical matrix $(b_{ij}=-b_{ji})$. This
proposition has the following motivation: we define the components of the
vector field {\bf V} by (10b) then, when we put this equation in (9), the
first term on the right hand side vanishes and we get the equation (10a).
This is not yet enough. We construct the curve as an integral curve of
the vector field with the components (10b), i.e., the solution of the
following set of equations (a non-autonomous system of differential
equations):
\begin{equation}
\frac{dx_i}{dt}=\sum\limits_{j=1}^{n-1}b_{ij}\frac{\partial E}{\partial
x_j}.
\end{equation}
Now with the solution of the equation (11) we have an explicit form of the
curve {\bf p}. And we know $E$, hence we know its partial derivatives.
Then for the function $f$ we can write down:
\begin{equation}
f=\int\left(\lim\limits_{\hat{\bf r}\rightarrow{\bf p}}\frac{\partial
E}{\partial t}\right)dt+const.
\end{equation}
So, with just the form of $E$ we can define the form of the other two
functions.

\underline{$\{{\bf p}\}$-case:} We just know the form of the curve.
However, for this case we require the following conditions:  $E\in
C^1(R^n,R)$, ${\bf p}\in C^2(R,R^n)$, $f\in C^1(R,R)$. We shall follow the
same methodology used in $\{E\}$-case. We know the explicit form of the
curve {\bf p}, hence we know its derivatives in an explicit way. We use
here a symbol $k_i(t)=dx_i/dt$ to denote these explicit functions. We have
the following two equations from the defining relation:
\begin{equation}
({\rm a})\;\frac{df}{dt}=
\lim\limits_{\hat{\bf r}\rightarrow{\bf
p}}\frac{\partial E}{\partial t}\qquad
{\rm and}\qquad({\rm b})\;\frac{\partial E}{\partial x_i}=
\sum\limits_{j=1}^{n-1}b_{ij}k_i(t),
\end{equation}
where $b_{ij}=const$ for all $i,j$. In this case we have supposed that the
components of the vector field, in the limit, are equal to the functions
$k_i(t)$. The solution to these equations is:
$$
E=\sum\limits_{i,j}^{n-1}b_{ij}k_i(t)x_i+T(t)\qquad{\rm and}\qquad
f=\int\lim\limits_{\hat{\bf r}\rightarrow{\bf
p}}\left\{\sum\limits_{i,j}^{n-1}b_{ij}\frac{dk_j}{dt}x_j\right\}dt
+\int\frac{dT}{dt}dt,
$$
where $T$ is an arbitrary function. In this case we have solved, first,
the equation (13b) and its solution $E$ is used to calculate the partial
derivative with respect to $t$. Then we have calculated the limit to get
the integrand to calculate $f$. Again, with just one entity, the curve, we
can define the other two functions in the functional  equation (6).

\underline{$\{f\}$-case:} We just know the form of $f$. The defining
relation is written as:
\begin{equation}
H(t)=\lim\limits_{\hat{\bf r}\rightarrow{\bf
p}}\left\{\sum\limits_{i=1}^{n-1}V_i(\hat{\bf r})\frac{\partial
E}{\partial x_i}+\frac{\partial E}{\partial t}\right\}.
\end{equation}
For this case we propose the following strategy (again we define the curve
as an integral curve of the vector field $V_i$):
\begin{equation}
({\rm a})\;(\forall i)\; \frac{dx_i}{dt}=H(t),\qquad ({\rm
b})\;H(t)\sum\limits_{i=1}^{n-1}\frac{\partial E}{\partial x_i}
+\frac{\partial E}{\partial t} =H(t).
\end{equation}
Hence the curve has the form $x_i=\int H(t)dt,\;(i=1,\ldots,n-1)$. The
function $E$ is determined by a first order partial differential equation
of a certain special form (Eq.(15b)).
\end{quotation}

One may think that the way in which we have solved the problems is
artificial because we introduced {\it ad hoc} vector fields in the
reasoning. This is not really the case, it is just the effect of our rigid
vision of the process of solution.

Consider, for example, the Poincar\'e-Cartan {\it 1-form} of classical
mechanics:
$$
W=\sum\limits_{i=1}^{n}p_idq_i-Hdt.
$$
We do not have the right to write down it as $W=dS$, where $S$ is the
action, until we prove that it is in fact an integrable {\it 1-form}. With
this purpose in mind we can attack the problem in the following way: we
suppose that the form is integrable and we write $p_i=\partial S/\partial
q_i$, $H(p_i,q_i)=-\partial S/\partial t$, and we get the Hamilton-Jacobi
equation. Hence, the problem of integrability is the problem of the
existence of solutions of the Hamilton-Jacobi equation. As it is
well-known, an analytic solution for this equation always exists locally
(Cauchy-Kovalevsky theorem), hence, the {\it 1-form} is a locally
integrable {\it 1-form}. In the dynamical problem we know the Hamiltonian
explicitly; but we know {\it neither} the form of the curve {\it nor} the
action as a function of the coordinates (not as a functional, because that
is another point of view). But, as it is well-known, if we can solve the
Hamilton-Jacobi equation we know the action and the solution of the
dynamical problem by means of a canonical transformation generated by this
action function. Clearly, in this case we have introduced all our
``auxiliary functions", the action and the Hamiltonian, to know the
explicit form of the curve in phase-space. Of course, we have required some
data: the form of the Hamiltonian and the supposition of integrability of
the {\it 1-form}. And from the theoretical point of view, it is enough to
construct the solution of the dynamical problem. However, we need to make
our distinction in this point: {\it the action as a function of the
coordinates differ from the function constructed by restriction of the
action to the curve}.

Another important point becomes clear when we use {\it 1-forms}: in all
the cases which we have treated, we need to suppose the integrability of a
{\it 1-form}. For example, when we treat the $\{E\}$-case we start from
the {\it 1-form}:
$$
dE=\sum\limits_i\frac{\partial E}{\partial x_i}dx_i+\frac{\partial
E}{\partial t}dt
$$
which is clearly integrable. Hence, we want to know a curve as an integral
curve of a vector field which we define as:
$$
X=\sum\limits_{i,j}b_{ij}\frac{\partial E}{\partial
x_j}\frac{\partial}{\partial x_i}+ \frac{\partial}{\partial t}.
$$
The inner product of these two tensors (the pairing between the tangent
and cotangent space) give us the result:
$$
\langle dE,X\rangle=\frac{\partial E}{\partial t}(x_1,\ldots,t),
$$
hence, the composition is in fact, the result of taking the limit of the
inner product in the integral curves of the vector field $X$. We can treat
the other cases from this point of view, but that is easy after this
explanation. In a geometric interpretation we have the following elements:
the tangent vectors, and the angle between them. In the $\{E\}$-case we
have the normal, but we have neither the tangent nor the angle; in the
$\{{\bf p}\}$-case the tangent, but we have neither the normal nor the
angle; finally, in the $\{f\}$-case we have the angle, but we have neither
the tangent nor the normal.

Now let us make a brief review of
the last case ($\{\}$-case).

The point is that in this case we have no any data and to treat it we need
some information. Heyting [3] notes that we ought to distinguish two
different concepts:
\begin{quotation}

1. Theories of the constructible.

2. Constructive theories

\end{quotation}
The first one is characterized by 3 conditions:
\begin{quotation}
\noindent
(a) we presuppose a mathematical theory in which the class of
constructible objects can be defined;\\
(b) the notion of a constructibility is no primitive;\\
(c) we have a liberty to choose the definition of a constructible, But, of
course, it must correspond to our intuitive notion of a mathematical
construction.
\end{quotation}
For the second point (the constructive theories) Heyting says: ``{\sl a
theory in which an object is only considered as existing after it has been
constructed. In other words, in a constructive theory there can be no
mentioning of other than constructible objects}". The main feeling of
Heyting is expressed in the following sentence: ``{\sl I am unable to give
an intelligible sense to the assertion that a mathematical object which has
not been constructed exists.}"

In the case which we want to treat we have no any data concerning the
equation $f=E\circ {\bf p}$. Hence, if we accept that we can only speak
about those objects which can be constructed explicitly (or, at least, we
have a method to construct them), the case which we are treating, the
\{\}-case, is not even a case. It is nothing, it is just a line of symbols
without any meaning. For this reason when one speaks about the functional
equation $f=E\circ {\bf p}$ one, in fact, is speaking about the
cases considered before: $\{E,{\bf p}\}$-, $\{E,f\}$-, $\{{\bf p},f\}$-,
$\{E\}$-, $\{f\}$-, and  $\{{\bf p}\}$-case.

As a last remark we can see that we have shown that the {\it generally
accepted} expressions of the type of Eq.(3) {\it cannot be valid}.

\section{About functional extensions}

We shall give our  problem the most general setting. Let us start with a
topological space $D$,  so that it is possible to construct the general
object of arrows:  $T(D,K)$ where $T$ is any covariant functor.  Hence we
can construct the functor:
\begin{equation}
T( D,\; ^*) : {\bf C}_1\rightarrow {\bf C}_2,
\end{equation}
where  ${\bf C}_i (i = 1, 2)$ are any small categories. Then for each
arrow we have $f\in T(D,K)$  the diagram: $f : D \rightarrow K$. For us,
the following situation is the most  important: the set $D$ is an object
with a given structure, so we use the symbol $P ( D )$
to denote its power set (which is a topology, of course, any topology is
a subset of the power set, but not any subset of the power set is a
topology). In this way, for each element in $P(D)$ we can define the
following elements:  $\langle f_A , A\rangle$  for all $A\in ( P( D )$.
Here the symbol $\langle f_A , A\rangle$ means that the object $A\in P(D)$
is put in correspondence with the function $f_A$. So
we may form the set:
\begin{equation} F_D = \{\langle f_A , A
\rangle|A\in P( D )\}
\end{equation}
of functional elements.

It  is clear that this procedure has been realized in a somewhat formal
manner, however, this is the more general form.  As we can see, there are
several elements which are important for our construction: the covariant
functor $T$,  the object $D$, its power set $P ( D )$, the set of
elements $F_D$ which is  the part in which the functions enter the
discussion and the two small categories: ${\bf C}_1, {\bf C}_2$.

\begin{quotation}
\noindent
{\small{\bf Definition 1:} We  shall call the symbol $\langle F_D,T,{\bf
C}_1, {\bf C}_2, P( D ) \rangle$ a general function.}
\end{quotation}

The idea behind a  general function is that all its elements are different
for each element of the power set of  $D$. Sometimes we can use a specific
topology instead of the power set, but this choice relies on our
convenience. Besides, we can see that in general, any topology is just a
subset of $P(D)$. The formation of a topology in the object $D$ can
respect its structure or not.  We use the categorical notions to introduce
the generality which they carry, because in general a function depends on
the categories in which it is defined, see [4] chap. 1, for more details.
Hence, the general setting is: {\it how are the different elements of a
general function related?} The problem may seem trivial without more
elaboration, however, as we have seen in the introduction,  in some realms
the problem is not trivial. Let us give a few additional examples.

{\bf Example 1:}  Consider the following  example (see  [5]), which is
clearly not trivial, suppose the following  choice:  $D = {\bf C}$ where
{\bf C} is the complex plane, if we use the symbol {\bf Anal}  to denote
the functor of the set of complex analytic functions we have:
\begin{equation}
{\bf Anal}({\bf C},\;^*) : {\bf Set}\rightarrow {\bf Set}
\end{equation}
or, to be more concrete, the arrows: $f : {\bf C \rightarrow C}$ of
complex analytic functions are at hand. We must consider the power set
$P({\bf C)}$ of the complex plane and the construction of the elements:
$\langle f_A , A \rangle$ for each element of the power set. This is the
most general situation for the choice that we have made of our basic
elements.  In this setting we have the following group of well-known
definitions [5]:
\begin{quotation}
\noindent
{\small Given two ``functional elements" $f(A) = \;_{\rm def}\langle f_A ,
A \rangle,\quad f(B) =\;_{\rm def}\langle f_B , B \rangle$,  we can say
that we have a direct analytic prolongation if, and only if, the following
two conditions hold:  $$ A\cap  B \neq \emptyset $$ $$ f_A = f_B \quad{\rm
in}\quad  A \cap B $$} \end{quotation}

So, we can  see that in general, the problem of analytic continuation is
a realization,  in the complex domain, of our definition of a general
function

{\bf Example 2:} Consider the functors:

\begin{equation}
     C(D,\;^* ) : {\bf Top\rightarrow CRng}\qquad {\rm and}\qquad
C^*(D,\;^*) : {\bf Top\rightarrow CRng}
\end{equation}
from  the topological spaces to the rings of continuous functions. The
functor  $C^*$ is for bounded functions. Here the problem is as follows: a
set $S$ is  $C$-embedded if, and only if, every function $f \in C( S )$
can be extended to a function ${\sl g} \in C( D )$. Here $S\subset  D$ and
$C(S)$ is an abbreviation of $C(S,S)$. The idea here is that the
extension is a $C$-function. The definition of $C^*$-embedding is similar.

One of the most important characteristics of the $C^*$-embedding is
Uryshon's {\bf lemma}:
\begin{quotation}
\noindent
{\small A subset  $S$ of the set $D$ is $C^*$-embedded in $D$ if, and only
if, any two completely separated sets in $S$ are completely separated in
$D$.}
\end{quotation}

We can see that  this lemma is just an assertion about functional
extensions, that  is, a theorem about the way in which the elements of a
general function  are related  [6], p. 18. Here the set $F_D$ can be
constructed once we  have fixed the topology of the spaces $S$ and $D$, or
at least the base of the  topology. If we use the power set we have a
conceptual generality, but we can fall into troubles for some purposes.
Let us take for topology of $D$ its power set $t(D)$, hence, the set $F_D$
can be formed and we have:  \begin{equation} F = \langle  F_D , C, {\bf
Top}, {\bf CRng} , t(D)\rangle \end{equation} as  our general function for
this case. Of course, we can construct $F$ without  recourse to the
Uryshon's lemma, however, this result gives us a way to  relate two
elements of the general function $F$.

{\bf Example 3:} Let  us come back to the example in the introduction.
Consider the functor:
\begin{equation}
C^{\infty} ( R^n ,\;^* ) : {\bf Vect\rightarrow Vect}
\end{equation}
so  the arrows: $f : R^n \rightarrow R^n$ where $R$ is the real line. The
power set is now $P( R^n )$,  and the set of functional elements is:
$\{\langle f_A , A \rangle| A \in P( R^n )\}$.  The notion of
differentiability does not change and we can define the derivative  of a
general function as the general function formed with the derivatives  of
the functional elements of the starting general function.  If one of such
elements is not differentiable, the general function is not.

\section{The limiting procedure}

Let us explain  in more detail the limiting procedure which can be used
for the elements  of a general function.  Consider the initial object $D$
and suppose a partition of the form:
\begin{equation}
D = \bigcup\limits_{i=1}^n G_i
\end{equation}
Hence,  the general function  is defined with the help of the elements of
the set $F_D = \{\langle f_i , G_i\rangle>\}$ and the functor $T( D ,\;^*
)$.  We can make this decomposition in many ways. For example: $D =
\bigcup_{A\in t(D)} A$, where $t( D )$ is the power set of $D$.

Now we define the system of sets:
\begin{equation}
P_i = \{ A\in P( D )| ( G_i \subset A  \}.
\end{equation}
In  other words: the set of all the sets  $A$ so that the set $G_i$ is
contained. Is very easy to show that each $P_i$ is a filter.
\begin{quotation}
\noindent
{\small{\bf Lemma:} Each $P_i$ is a model of a  filter in  $D$  (see [6]
p.  24).}
\end{quotation}

{\tt Proof}: $\Diamond$ (a) We  can see that $\emptyset$ is not an element
of $P_i$, any $i$, because if $\emptyset \in P_i$ then we can find a set
$A$ such that $A\subset \emptyset$ which is a contradiction, hence we have
proved the first axiom.  (b) If we suppose that $A, B \in P_i$ then $A\cap
B \in P_i$ because, at least $A$ and $B$ have the set $G_i$ in common,
hence $G_i$ is in their intersection, but this is the condition for
belonging to $P_i$.  The second axiom is satisfied. (c) If $A\in P_i$, $B
\subset D$ and $A\subset B$ it is very easy to see that $B\in P_i$. The
lemma is proved.$\Diamond$

This lemma  (such trivial as it is) is important, because with a filter we
 can define a limit for the elements of a general function.  In fact,
 given the filter $P_i$ of the element $i$ in the partition, we have a
function $f_i$ which maps the element $G_i$. Clearly, the elements $f_i (
G_i )$ are the images of the set $D$ by the general function $F_D$. So,
we define the filter of the image of $D$ by the general function $F_D$ as
$F_D ( Pi ) = \{ A\in K | f_i ( Gi ) \subset A \}$. This is clearly a
filter.  With these elements it is possible to set up a well-known
definition of the limiting procedure for the elements of a general
function.

\begin{quotation}
\noindent
{\small{\bf Definition 2:} A set  $G_i$ is the limit of a filter $H$ if,
and only if, $H$ is {\it stronger} than the filter $P_i$.}
\end{quotation}
\begin{quotation}
\noindent
{\small{\bf Definition 3:} Consider  the filter $P_i$,  hence the set
$A\subset K$ is the limit of the general function $F_D$  under the filter
$P_i$ if, and only if the set $A$ is the limit of the filter $F_D( P_i )$.
That is, without abbreviations:  the filter $F_D( P_i )$ is stronger than
the filter formed with the sets that contain $A$.}
\end{quotation}

We say that a  filter $H$ is stronger than the filter $B$ (of course, both
filters defined on the same space) if, and only if, for any $a\in B$ there
is a set $b\in H$ such that $b\subset a$. Of course, this is just the
notion of approximation, because a filter $H$ is stronger than a filter
$B$ if their elements are nearest to a certain set than the elements of
$B$.  Now let us use this concept for the example in the introduction. We
have the equation:
\begin{equation}
\frac{d}{dt} E_{\gamma} = \lim\limits_{{\bf r}\rightarrow {\bf r}(t)}
\left\{({\bf V}\cdot\nabla)E_A + \frac{\partial}{\partial t}E_A\right\},
\end{equation}
where  $E_{\gamma}$ is a function along the curve $\gamma$ and $E_A$ is a
function defined on the set  $A$. Now, let us give a precise meaning to
the process involved.  We have the  general function $E_D$  and two of
their functional elements are involved:  $\langle
E_{\gamma},\gamma\rangle$, $\langle E_A, A \rangle$ where $\gamma$ and
$A$ are sets in $D$.  Hence we can see that the limiting procedure affects
only the functional element $\langle ({\bf V}\cdot\nabla)E_A , A \rangle$
so, we do the following:  we select a set $d \subset \gamma$ and we form
its filter $P_d$ , so, a set $B \in R^n$ in the image of the functional
element $({\bf V}\cdot\nabla)E_A$ is its limit if, and only if the
filter formed with the sets that contain the image of  $({\bf
V}\cdot\nabla)E_A$  is stronger than the filter formed with the sets that
contain $B$.  Of course the extension of this definition covers the usual
$\varepsilon$-$\delta$ arguments.

\section{Conclusions}

As  promised in  the introduction, we have solved in its most general form
the  ambiguity  which arises in the differential calculus of several
variables with the help of  category  theory. Besides we have showed
several examples of realizations of our  construction.

\acknowledgments

The authors would like to express their gratitude to Prof. Valeri
Dvoeglazov for his discussions and critical comments. We would also like
to thank Annamaria D'Amore for revising the manuscript.

 

\begin{references}

\bibitem{1} V. Y. Arnold, {\it Matematicheskie Metody Classicheskoi
Mehaniki} (Nauka, Moscow, 1989) [Englishish translation:  {\it
Mathematical Methods of Classical Mechanics} (Springer-Verlag, New York,
1989)].



\bibitem{2} A. E. Chubykalo and  R. Smirnov-Rueda, Modern Physics
Letters A {\bf 12}(1), 1 (1997).

\bibitem{3} A. Heyting, {\it Constructivity in Mathematics}
(North-Holland, Amsterdam, 1956).


\bibitem{4} S. MacLane, {\it Categories for the Working
Mathematician}  (Springer-Verlag, New York, 1971).



\bibitem{5} L. Ahlfors, {\it Complex Analysis} (Mc Graw-Hill, New
York, 1953)



\bibitem{6} L. Gillman and M. Jerison, {\it Rings of continuous functions}
(D. Van Nostrand Company, New York,  1960).





\end{references}
\end{document}